\documentclass[12pt]{amsart}
\usepackage{amsfonts}
\usepackage{amssymb}
\usepackage{amscd}

\newcommand{\Eqref}[1]{(\ref{#1})}

\newlength\Itemwidth
\newlength\Itemmargin

\newenvironment{List}[1]%
{\settowidth{\Itemwidth}{#1}%
\begin{list}{}%
{%
\setlength{\labelsep}{6pt}%
\setlength{\itemsep}{1pt}%
\setlength{\parsep}{0pt}%
\setlength{\itemindent}{0pt}%
\setlength{\labelwidth}{\Itemwidth}%
\setlength{\leftmargin}{\Itemwidth}%
\addtolength{\leftmargin}{\labelsep}%
}}{\end{list}}

\newtheorem{theorem}[subsection]{Theorem}

\newtheorem{lemma}[subsection]{Lemma}

\newtheorem{PreExamples}[subsection]{Examples}
\newtheorem{PreExample}[subsection]{Example}
\newtheorem{PreRemarks}[subsection]{Remarks}
\newtheorem{PreRemark}[subsection]{Remark}
\newenvironment{examples}%
{\begin{PreExamples}\rm}{\end{PreExamples}}
\newenvironment{example}%
{\begin{PreExample}\rm}{\end{PreExample}}
\newenvironment{remarks}%
{\begin{PreRemarks}\rm}{\end{PreRemarks}}
{\begin{PreRemark}\rm}{\end{PreRemark}}

\newtheorem{PreState}[subsection]{\statetype}
\newenvironment{state}[1]%
{\def\statetype{#1} \begin{PreState}\rm}{\end{PreState}}

\newcommand{\rank}{\mathop{\mathrm{rk}}\nolimits}
\newcommand{\Card}{\sharp}
\newcommand{\highlight}[1]{{\sl #1}}

\begin{document}

\title[Sum of the index]{On the sum of the index of a parabolic subalgebra and of 
its nilpotent radical}
\author{Rupert W.T. Yu}
\dedicatory{UMR 6086 du C.N.R.S.,
D\'epartement de Math\'ematiques\\
Universit\'e de Poitiers, 
T\'el\'eport 2 -- BP 30179,
Boulevard Marie et Pierre Curie,\\
86962 Futuroscope Chasseneuil Cedex,
France.}
\address{UMR 6086 du C.N.R.S.\\
D\'epartement de Math\'ematiques\\
Universit\'e de Poitiers\\
T\'el\'eport 2 -- BP 30179\\
Boulevard Marie et Pierre Curie\\
86962 Futuroscope Chasseneuil Cedex\\
France.}
\email{yuyu@math.univ-poitiers.fr}

\begin{abstract}
In this short note, we investigate the following question of Panyushev stated in \cite{Pany}:
``Is the sum of the index of a parabolic subalgebra of a semisimple Lie algebra $\mathfrak{g}$
and the index of its nilpotent radical always greater than or equal to the rank of 
$\mathfrak{g}$?''. Using the formula for the index of parabolic subalgebras
conjectured in \cite{TYseaweed} and proved in \cite{MFJ,Jindex}, we give a positive 
answer to this question. Moreover, we also obtain a necessary and 
sufficient condition for this sum to be equal to the rank of $\mathfrak{g}$. This provides
new examples of direct sum decomposition of a semisimple Lie algebra 
verifying the ``index additivity condition'' as stated by Ra\"\i s in \cite{Rais}.
\end{abstract}

\maketitle

\section{Introduction}

Let $\mathfrak{g}$ be a Lie algebra over an algebraically closed field
$\Bbbk$ of characteristic zero. For $f \in \mathfrak{g}^{*}$, we denote by
$\mathfrak{g}^{f}=\{ X\in \mathfrak{g} ; f([X,Y])= 0$ for all $Y \in \mathfrak{g} \}$,
the annihilator of $f$ for the coadjoint representation of $\mathfrak{g}$. The
\highlight{index} of $\mathfrak{g}$, denoted by $\chi ( \mathfrak{g})$, is defined 
to be 
$$
\chi (\mathfrak{g})=\min_{f\in \mathfrak{g}^{*}} \dim \mathfrak{g}^{f}.
$$
It is well-known that if $\mathfrak{g}$ is an algebraic Lie algebra and $G$ denote
its algebraic adjoint group, then $\chi (\mathfrak{g})$ is the transcendence degree
of the field of $G$-invariant rational functions on $\mathfrak{g}^{*}$.

The index of a semisimple Lie algebra $\mathfrak{g}$ is equal to the rank of $\mathfrak{g}$.
This can be obtained easily from the isomorphism between $\mathfrak{g}$ and $\mathfrak{g}^{*}$
via the Killing form. There has been quite a lot of recent work on the determination 
of the index of certain subalgebras of a semisimple Lie algebra: 
parabolic subalgebras and related subalgebras 
(\cite{DK}, \cite{Panyseaweed}, \cite{TYseaweed}, \cite{Mor2}), 
centralizers of elements and related subalgebras 
(\cite{Pany}, \cite{Char}, \cite{Yaki}, \cite{Mor1}).

Let $\mathfrak{g}$ be a semisimple Lie algebra, $\mathfrak{p}$ a parabolic subalgebra of 
$\mathfrak{g}$ and $\mathfrak{u}$ (resp. $\mathfrak{l}$) the nilpotent radical 
(resp. a Levi factor) of $\mathfrak{p}$. 
In \cite[Corollary 1.5 (i)]{Pany}, Panyushev showed that
\begin{equation}\label{oldupperbound}
\chi (\mathfrak{p})+\chi (\mathfrak{u}) \leq \dim \mathfrak{l}.
\end{equation}
He then suggested \cite[Remark (ii) of Section 6]{Pany} that 
\begin{equation}\label{indexsumformula}
\chi (\mathfrak{p})+\chi (\mathfrak{u}) \geq \rank \mathfrak{g}.
\end{equation}
For example, it is well-known that if $\mathfrak{b}$ is a Borel subalgebra of $\mathfrak{g}$
and $\mathfrak{n}$ is its nilpotent radical, then 
$\chi (\mathfrak{b})+\chi (\mathfrak{n})=\rank \mathfrak{g}$ (see for example
\cite{TYstable}, \cite[Chapter 40]{TYbook}). It is therefore also interesting to characterise
parabolic subalgebras where equality holds in \Eqref{indexsumformula}. Indeed, in
\cite{Rais}, Ra\"\i s looked for examples of direct sum decompositions 
$\mathfrak{g}=\mathfrak{m}\oplus \mathfrak{n}$ verifying the ``index
additivity condition'', namely $\mathfrak{m}$ and $\mathfrak{n}$ are Lie subalgebras of 
$\mathfrak{g}$ and 
$$
\chi (\mathfrak{g})=\chi (\mathfrak{m}) + \chi (\mathfrak{n}).
$$

If $\mathfrak{u}_{-}$ denotes the nilpotent radical of the opposite parabolic subalgebra
$\mathfrak{p}_{-}$ of $\mathfrak{p}$, then $\mathfrak{g}=\mathfrak{p}\oplus \mathfrak{u}_{-}$
and the Lie algebras $\mathfrak{u}$ and $\mathfrak{u}_{-}$ are isomorphic. Thus 
parabolic subalgebras such that equality holds  in \Eqref{indexsumformula}
would provide examples of direct sum decompositions verifying the index
additivity condition.

Using the formula, conjectured in \cite{TYseaweed} and proved in
\cite{MFJ,Jindex}, for the index of parabolic subalgebras, 
we obtain a formula for the sum $\chi (\mathfrak{p}) + \chi (\mathfrak{u})$. 
By a careful analysis of root systems, we prove the inequality \Eqref{indexsumformula} and give
a necessary and sufficient condition of equality to hold in \Eqref{indexsumformula} (See Theorem
\ref{Main}).

To describe the index of a parabolic subalgebra, and the index of its nilpotent radical
we need to recall Kostant's cascade construction of pairwise strongly
orthogonal roots (\cite{Jan}, \cite{Jos}, \cite{TYbook}). 

Let us fix a Cartan subalgebra $\mathfrak{h}$ of $\mathfrak{g}$ and a
Borel subalgebra $\mathfrak{b}$ of $\mathfrak{g}$ containing $\mathfrak{h}$. Denote 
by $R$, $R^{+}$ and $\Pi=\{ \alpha_{1},\dots ,\alpha_{\ell} \}$ 
respectively the set of roots, positive roots and simple roots with 
respect to $\mathfrak{h}$ and $\mathfrak{b}$. For any $\alpha \in R$, let $\mathfrak{g}_{\alpha}$
be the root subspace associated to $\alpha$. Choose $X_{\alpha}$ such that 
$\alpha ([X_{\alpha},X_{-\alpha}])=2$. We shall write 
$\alpha^{\vee}=[X_{\alpha}, X_{-\alpha}]\in \mathfrak{h}$,
and for $\lambda \in \mathfrak{h}^{*}$, $\langle \lambda, \alpha^{\vee} \rangle =\lambda
(\alpha^{\vee})$. For $S\subset \Pi$, we denote by
$R_{S}=R\cap \mathbb{Z}S$, $R_{S}^{+}=R_{S} \cap R^{+}$. If $S$ is connected,
then we shall denote by $\varepsilon_{S}$ the highest root of $R_{S}$.

Let $S\subset \Pi$. We define $\mathcal{K}(S)$ inductively as follows:
\begin{List}{(iii)}
\item[a)] $\mathcal{K}(\emptyset)=\emptyset$.
\item[b)] If $S_{1},\dots ,S_{r}$ are the connected components of $S$, then
$
\mathcal{K}(S)=\mathcal{K}(S_{1})\cup \cdots \cup \mathcal{K}(S_{r}).
$
\item[c)] If $S$ is connected, then 
$
\mathcal{K}(S)=\{ S \} \cup \mathcal{K}(\widehat{S})
$
where $\widehat{S}=
\{ \alpha \in S ; \langle \alpha , \varepsilon_{S}^{\vee} \rangle=0 \}$.
\end{List}

It is well-known that (see for example \cite[Chapter 40]{TYbook}) elements of $\mathcal{K}(S)$
are connected subsets of $S$. Moreover, if we denote by 
$\mathcal{R}(S)=\{ \varepsilon_{K} ; K\in \mathcal{K}(S) \}$, then $\mathcal{R}(S)$
is a maximal set of pairwise strongly orthogonal roots in $R_{S}$. 

Let us also recall the following properties of $\mathcal{K}(S)$:

\begin{lemma}\label{lemma}
Let $S$ be a subset of $\Pi$, $K,K'\in \mathcal{K}(S)$ and set 
$$
\begin{array}{rl}
\Gamma^{K}& =\{ \alpha \in R_{K} ; \langle \alpha ,\varepsilon_{K}^{\vee} \rangle > 0 \}\\
& =\{ \alpha= \sum_{\beta\in K} n_{\beta}\beta \in R_{K}^{+} ; n_{\beta} > 0 \hbox{ for some }
\beta\in K\setminus \widehat{K} \}.
\end{array}
$$
\begin{List}{(iii)}
\item[i)] We have either $K\subset K'$ or $K'\subset K$ or $K$ and $K'$ are connected components
of $K\cup K'$.
\item[ii)] $\Gamma^{K} = R_{K}^{+} \setminus \{ \beta \in R_{K}^{+} ; 
\langle \beta , \varepsilon_{K}^{\vee} \rangle =0 \}$. In particular, 
$R_{K}^{+}$ is the disjoint union of the $\Gamma^{K}$'s, $K\in \mathcal{K}(S)$.
\item[iii)] $\sum_{\alpha \in \Gamma^{K}} \mathfrak{g}_{\alpha}$ 
is a Heisenberg Lie algebra whose centre is
$\mathfrak{g}_{\varepsilon_{K}}$. Thus if $\alpha ,\beta \in \Gamma^{K}$ verify
$\alpha + \beta \in R$, then $\alpha +\beta =\varepsilon_{K}$.
\item[iv)] Suppose that $\alpha \in \Gamma^{K}$ and $\beta \in \Gamma^{K'}$ verify
$\alpha + \beta \in R$, then either $K\subset K'$ and $\alpha + \beta \in \Gamma^{K'}$
or $K'\subset K$ and $\alpha + \beta \in \Gamma^{K}$.
\end{List}
\end{lemma}

\goodbreak
\begin{state}{Table}\label{table}
The cardinality of $\mathcal{K}(\Pi )$ is listed below for an irreducible root system $R$:
\begin{center}
\begin{tabular}{|c|c|c|c|c|c|c|c|c|c|}\hline
Type & $A_{\ell}, \ell \geq 1$ &
$B_{\ell}, C_{\ell}, \ell \geq 2$ &
$D_{\ell} , \ell \geq 4$ & $E_{6}$ & $E_{7}$ & $E_{8}$ & $F_{4}$
& $G_{2}$ \\ \hline
$\Card \mathcal{K}(\Pi )$ &
\strut\rule[-2.8ex]{0ex}{7ex}
$\left[\displaystyle\frac{\ell+1}{2}\right]$ & $\ell$ & $2 \left[
\displaystyle\frac{\ell}{2}\right]$& $4$ & $7$ & $8$ & $4$ & $2$ \\ \hline
\end{tabular}
\end{center}
where for any $x\in\mathbb{Q}$, $[x]$ is the unique integer such that
$[x]\leq x<[x]+1$.
\end{state}

\begin{examples}\label{examples}
Let $R$ be an irreducible root system. We shall use the numbering of simple roots in 
\cite[Chapter 18]{TYbook}. Set $k=\Card \mathcal{K}(\Pi )$. 

1) Let $R$ be of type $A_{\ell}$. Then
$\mathcal{K}(\Pi )=\bigcup_{i=1}^{k}\{ K_{i}\}$ where
$K_{i}=\{ \alpha_{i},\dots ,\alpha_{\ell +1-i} \}$. For $1\leq i\leq k$,
$$
\Gamma^{K_{i}} = \{ \alpha_{i}+ \cdots + \alpha_{i+r}\ , \ 
\alpha_{\ell+1-i-r} + \cdots + \alpha_{\ell+1-i}
; 0\leq r \leq \ell -2i \}  \cup \{ \varepsilon_{K_{i}} \}.
$$

2) Let $R$ be of type $D_{2n+1}$. Then $k=2n$, 
$$
\mathcal{K}(\Pi )=\bigcup_{i=1}^{n} \{ K_{i}\} \cup \{ L_{i}\}
$$
where 
$K_{i}=\{ \alpha_{2i-1}, \dots ,\alpha_{2n+1} \}$
and 
$L_{i}=\{ \alpha_{2i-1} \}$. For $1\leq i\leq n$, 
$$
\Gamma^{K_{i}} = \bigg\{ \sum_{j=2i-1}^{\ell} m_{j}\alpha_{j} ; m_{2i}\neq 0\bigg\} \ , \ 
\Gamma^{L_{i}} = \{ \alpha_{2i-1} \}.
$$

3) Let $R$ be of type $E_{6}$. Then
$$
\mathcal{K}(\Pi )=\{ \Pi \} \cup 
\{ \{ \alpha_{1},\alpha_{3},\alpha_{4},\alpha_{5},\alpha_{6} \} \}
\cup 
\{ \{ \alpha_{3} , \alpha_{4},\alpha_{5} \} \} 
\cup 
\{ \{ \alpha_{4} \} \} .
$$
\end{examples}

\section{Main result}\label{main}

Recall that for any subset $S\in\Pi$,
$$
\mathfrak{p}_{S}=\mathfrak{h}\oplus \bigoplus_{\alpha\in R_{S}\cup R^{+}} \mathfrak{g}_{\alpha}
$$
is a (standard) parabolic subalgebra of $\mathfrak{g}$.
Any parabolic subalgebra of $\mathfrak{g}$ is conjugated to a standard
parabolic subalgebra. The Lie subalgebra
$$
\mathfrak{u}_{S}=\bigoplus_{\alpha\in R^{+}\setminus R_{S}} \mathfrak{g}_{\alpha}
$$
is the nilpotent radical of $\mathfrak{p}_{S}$.

For $S\subset \Pi$, 
denote by 
$V_{S}$ the vector subspace of $\mathfrak{h}^{*}$ spanned by the elements
of $\mathcal{R}(S)$ and $\mathcal{R}(\Pi )$. Set
$$
\begin{array}{c}
\mathcal{E}_{S}=\{ K \in \mathcal{K}(\Pi ) ; X_{\varepsilon_{K}} \in \mathfrak{u}_{S} \}
=\{ K \in \mathcal{K}(\Pi ) ; \varepsilon_{K} \not\in R_{S}\} \\
\hbox{ and } \
\mathcal{Q}_{S}=\bigg( \bigcup_{K\in \mathcal{E}_{S}} \Gamma^{K}\bigg) \cap R_{S}^{+}.
\end{array}
$$
The subset $\mathcal{E}_{S}$ has the following simple characterisation.

\begin{lemma}\label{TSlemma}
Let $T_{S}$ be the union of  $K\in \mathcal{K}(\Pi )$ verifying $K\subset S$. Then
$\mathcal{E}_{S}=\mathcal{K}(\Pi ) \setminus \mathcal{K}(T_{S})$.  
\end{lemma}
\begin{proof}
This is straightforward.
\end{proof}

Our main result is the following theorem.

\begin{theorem}\label{Main}
Let $S\subset \Pi$. Then
\begin{equation}\label{newupperbound}
\chi (\mathfrak{p}_{S})+\chi (\mathfrak{u}_{S})
=\rank \mathfrak{g} +  \Card \mathcal{K}(S) - \Card \mathcal{K}(T_{S}) +
2 (\Card \mathcal{K}(\Pi )-\dim V_{S} ) + \Card \mathcal{Q}_{S}.
\end{equation}
We have $\chi (\mathfrak{p}_{S})+\chi (\mathfrak{u}_{S})\geq \rank \mathfrak{g}$,
and equality holds
if and only if the following conditions are satisfied:
\begin{List}{(iii)}
\item[i)] $\Card ( \mathcal{K}(S)\cup \mathcal{K}(\Pi ) )= \dim V_{S}$.
\item[ii)] For any connected component $S'$ of $S$, 
we have either $S'\in \mathcal{K}(\Pi)$ or $\Card (S'\setminus T_{S})=1$.
\end{List}
\end{theorem}
\begin{proof} 
The formula for the sum $\chi (\mathfrak{p}_{S})+\chi (\mathfrak{u}_{S})$
is a direct consequence of the formula of the index of parabolic subalgebras
conjectured in \cite{TYseaweed} and proved
in \cite{MFJ,Jindex} :
\begin{equation}\label{paraformula}
\chi (\mathfrak{p}_{S}) = \rank \mathfrak{g} + \Card \mathcal{K}(\Pi) +
\Card \mathcal{K}(S) - 2 \dim V_{S} \ , 
\end{equation}
and the formula for the index of $\mathfrak{u}_{S}$ 
(see for example \cite[Chapter 40]{TYbook}), which, in view of Lemma \ref{lemma}, 
can be expressed in the following way: 
\begin{equation}\label{radformule}
\chi (\mathfrak{u}_{S})=
\Card \mathcal{E}_{S} + \sum_{K \in \mathcal{E}_{S}} \Card \Gamma^{K} - \dim \mathfrak{u}_{S}
= \Card \mathcal{E}_{S} +  \Card \mathcal{Q}_{S}
\end{equation}

To prove the rest of the theorem, we may clearly 
assume that $\mathfrak{g}$ is simple. Observe that 
\begin{equation}\label{S-T}
\Card \mathcal{K}(S) - \Card \mathcal{K}(T_{S}) \geq 0.
\end{equation}
\vskip0.5em

1) Let $S_{1},\dots ,S_{r}$ be the connected components of $S$.
For each $i$, there is a unique $K_{i}\in \mathcal{K}(\Pi )$ (see Lemma \ref{lemma})
such that $\varepsilon_{S_{i}} \in \Gamma^{K_{i}}$. 
If $K_{i}=S_{i}$, then $S_{i}$ is a connected component of $T_{S}$. 
Otherwise $K_{i} \in \mathcal{E}_{S}$, and we have
$$
\varepsilon_{S_{i}} \in \Gamma^{K_{i}} \cap \Gamma^{S_{i}} \subset \mathcal{Q}_{S}.
$$

2) It follows from Point 1) that $\mathcal{Q}_{S}=\emptyset$ if and only if
$\mathcal{K}(S)\subset \mathcal{K}(\Pi )$ (or equivalently $S=T_{S}$).
\vskip0.5em

3) Note that the connected components of $T_{S}$ are the connected components of 
$T_{S}\cap S_{i}$, it follows again from Point 1) that:
$$
\begin{array}{rcl}
\Card \mathcal{K}(S)-\Card \mathcal{K}(T_{S}) 
& = & \displaystyle\sum_{i=1}^{r} (\Card \mathcal{K}(S_{i}) -\Card \mathcal{K}(T_{S}\cap S_{i}) )\\
& = & \displaystyle\sum_{K_{i} \in \mathcal{E}_{S}} 
( \Card \mathcal{K}(S_{i}) - \Card \mathcal{K}(T_{S}\cap S_{i}) ). 
\end{array}
$$

4) From Table \ref{table}, we have $\dim V_{S} = \Card \mathcal{K}(\Pi )=\rank \mathfrak{g}$ 
in the cases where $\mathfrak{g}$ is of type
$B_{\ell}$, $C_{\ell}$, $D_{2n}$, $E_{7}$, $E_{8}$, $F_{4}$ et $G_{2}$. 
The inequality follows immediately from \Eqref{newupperbound} and \Eqref{S-T}, and
the condition for equality follows from Point 2). 
\vskip0.5em

5) \underbar{Type $A_{\ell}$}. 

For any $i$ verifying $S_{i}\neq K_{i}$,  by Point 1), Lemma \ref{lemma} 
and Examples \ref{examples}, half of $\Gamma^{S_{i}}\setminus \{  \varepsilon_{S_{i}} \}$ belongs
to $\mathcal{Q}_{S}$. Since $\Card (\Gamma^{S_{i}})=2\Card (S_{i})-1$ (Examples \ref{examples}), 
such an $S_{i}$ contributes $\Card (S_{i})$ elements of $\mathcal{Q}_{S}$.

Again, since we are in type $A$, $\mathcal{K}(\Pi )$ is totally ordered by inclusion. It
follows that $T_{S}$ is connected. Without loss of generality, we may assume that
$T_{S}\subset S_{1}$. 

Suppose that $S_{1}= T_{S}$. Then from the previous discussion, we deduce that 
$$
\chi (\mathfrak{p}_{S})+\chi (\mathfrak{u}_{S}) \geq 
\rank \mathfrak{g} + \Card \mathcal{K}(S\setminus S_{1})
+ 2 (\Card \mathcal{K}(\Pi )- \dim V_{S}) + \Card (S\setminus S_{1}).
$$
But our hypothesis implies that 
\begin{equation}\label{eqnA}
\dim V_{S} \leq \Card \mathcal{K}(\Pi ) + \Card \mathcal{K}(S \setminus S_{1})
=\Card (\mathcal{K}(\Pi )\cup \mathcal{K}(S) ),
\end{equation}
so
$$
\chi (\mathfrak{p}_{S})+\chi (\mathfrak{u}_{S}) 
\geq 
\rank \mathfrak{g} + \Card (S \setminus S_{1} ) -  \Card \mathcal{K}(S\setminus S_{1}).
$$
Hence $\chi (\mathfrak{p}_{S})+\chi (\mathfrak{u}_{S}) \geq \rank \mathfrak{g}$. 
For equality to hold, we must have equality in \Eqref{eqnA} and
$$
\Card (S \setminus S_{1} )=\Card \mathcal{K}(S\setminus S_{1}).
$$
This latter is only possible if $\Card (S_{i})=1$ for $i\geq 2$, so we have conditions 
(i) and (ii). Conversely, suppose that conditions (i) and (ii) are verified, then 
$\Card (S_{i})=1$ for $i\geq 2$. Consequently $\mathcal{Q}_{S}=S\setminus S_{i}$ by Point 1)
and the definition of $\mathcal{Q}_{S}$. 

Suppose that $S_{1}\supsetneq T_{S}$ (this includes the case $T_{S}=\emptyset$). Then
$\mathcal{K}(S)\cap \mathcal{K}(\Pi )=\emptyset$.
Thus 
\begin{equation}\label{eqnA2}
\dim V_{S} \leq \Card \mathcal{K}(\Pi ) + \Card \mathcal{K}(S)
=\Card ( \mathcal{K}(\Pi )\cup \mathcal{K}(S) ).
\end{equation} 
We deduce from Point 1) and 
the remark in the first paragraph of Point 5) that
$$
\begin{array}{rl}
\chi (\mathfrak{p}_{S})+\chi (\mathfrak{u}_{S}) & \geq  
\rank \mathfrak{g} + \Card \mathcal{K}(S) - \Card \mathcal{K}(T_{S})
+ 2 (\Card \mathcal{K}(\Pi )- \dim V_{S}) + \Card S \\
& \geq \rank \mathfrak{g} + \Card (S) - \Card \mathcal{K}(S) - \Card \mathcal{K}(T_{S}) .
\end{array}
$$
Hence
$$
\chi (\mathfrak{p}_{S})+\chi (\mathfrak{u}_{S}) 
\geq \rank \mathfrak{g} +  \displaystyle \sum_{i=1}^{r} \left( 
\Card (S_{i}) - \left[ \frac{\Card (S_{i}) +1}{2} \right] \right) - 
\displaystyle\left[ \frac{\Card (T_{S}) +1}{2} \right] .
$$
Since $T_{S}\subsetneq S_{1}$, we deduce from Table \ref{table} that
$$
\Card (S_{1}) - \displaystyle\left[ \frac{\Card (S_{1}) +1}{2} \right]  - 
\displaystyle\left[ \frac{\Card (T_{S}) +1}{2} \right] \geq 0.
$$
So we have our inequality $\chi (\mathfrak{p}_{S})+\chi (\mathfrak{u}_{S}) \geq \rank \mathfrak{g}$.

Now for the equality $\chi (\mathfrak{p}_{S})+\chi (\mathfrak{u}_{S}) =\rank \mathfrak{g}$
to hold, we must have equality in \Eqref{eqnA2} and
$\Card \mathcal{Q}_{S}=\Card (S)$, 
$$
\Card (S_{1}) - \left[ \frac{\Card (S_{1}) +1}{2} \right] - 
\displaystyle\left[ \frac{\Card (T_{S}) +1}{2} \right]=0
\hbox{ and }
\Card (S_{i}) - \left[ \frac{\Card (S_{i}) +1}{2} \right]=0
$$ 
for $i\geq 2$. This implies that $\Card (S_{1})=\Card (T_{S})+1$, and $\Card (S_{i})=1$ for
$i\geq 2$. So we have conditions (i) and (ii). Conversely, 
if conditions (i) and (ii) are verified,
then $\Card (S_{i})=1$ for $i\geq 2$, and $\Card (S_{1})=\Card (T_{S})+1$. In view of the above,
to show that $\chi (\mathfrak{p}_{S})+\chi (\mathfrak{u}_{S}) =\rank \mathfrak{g}$, it suffices to check
that $\Card (\mathcal{Q}_{S}\cap R_{S_{1}}^{+})=\Card (S_{1})$, which is a straightforward
verification. 
\vskip0.5em

6) \underbar{Type $D_{2n+1}$}.

In this case, $\Card \mathcal{K}(\Pi )=\rank \mathfrak{g}-1$. Let us use the numbering of 
simple roots in \cite[Chapter 18]{TYbook}. We check easily that 
$\alpha_{1},\dots ,\alpha_{\ell-2}\in V_{\Pi}$. 

If $\dim V_{S}=\Card \mathcal{K}(\Pi )$, then the inequality follows from 
\Eqref{newupperbound} and \Eqref{S-T}, and the condition for equality follows from Point 2).

Suppose now that $\dim V_{S} = \rank \mathfrak{g}$ and
$\alpha_{\ell-1}\in S$ (the case $\alpha_{\ell}\in S$ being analogue).
Then
\begin{equation}\label{eqnD1}
\chi (\mathfrak{p}_{S})+\chi (\mathfrak{u}_{S})  
= \rank \mathfrak{g} + \Card \mathcal{K}(S) - \Card \mathcal{K}(T) -
2 + \Card (\mathcal{Q}_{S}),
\end{equation}
and the connected component $S_{1}$ of $S$ containing $\alpha_{\ell-1}$ is not in 
$\mathcal{K}(\Pi )$ (for otherwise, we would have
$\dim V_{S} = \Card \mathcal{K}(\Pi )$).

By Point 1), $\varepsilon_{S_{1}} \in \mathcal{Q}_{S}$. By examining
the possibilities for $S_{1}$ and $K_{1}$ (Examples \ref{examples}), we verify that
\begin{equation}\label{eqnD2}
\Card \mathcal{K}(S_{1}) - \Card \mathcal{K}(T_{S}\cap S_{1}) +
\Card (\Gamma^{K_{1}} \cap \Gamma^{S_{1}}) \geq 2
\end{equation}
with equality if and only if 
$S_{1}$ is of type $A_{1}$ or $A_{2}$. We have therefore obtained the inequality.

In fact, we showed in the previous paragraph that already we have
$$
\rank \mathfrak{g} + \Card \mathcal{K}(S_{1}) - \Card \mathcal{K}(T\cap S_{1}) - 2 +
\Card (\Gamma^{K_{1}} \cap \Gamma^{S_{1}}) \geq \rank \mathfrak{g}.
$$
So if $\chi (\mathfrak{p}_{S})+\chi (\mathfrak{u}_{S}) =\rank \mathfrak{g}$, 
then from \Eqref{eqnD1} and the above inequality, 
we must have $\mathcal{K}(S\setminus S_{1}) \subset \mathcal{K}(\Pi )$,
and also equality in \Eqref{eqnD2}. Hence conditions (i) and (ii).
Conversely, suppose that conditions (i) and (ii) are verified, then
the fact that $\alpha_{1},\dots ,\alpha_{\ell-2}\in V_{\Pi}$ implies that
$\mathcal{K}(S\setminus S_{1}) \subset \mathcal{K}(\Pi )$ and 
$\Card \mathcal{K}(S_{1})=1$. Hence $S_{1}$ is of type $A_{1}$, $A_{2}$.
It is then easy to check that 
$\chi (\mathfrak{p}_{S})+\chi (\mathfrak{u}_{S}) =\rank \mathfrak{g}$.
\vskip0.5em

7) \underbar{Type $E_{6}$}.

Here, we have $\Card \mathcal{K}(\Pi )=4$ and $\alpha_{2},\alpha_{4} \in V_{\Pi}$.
Let $S_{1}$ be a connected component of  $S$ such that
$\dim V_{S_{1}} > 4$. Under these conditions, the possibilities are:
$$
\begin{array}{|c|c|c|c|c|c|c|c|}
\hline 
S_{1} & \dim V_{S_{1}} & T_{S}\cap S_{1} & 
\Card \mathcal{K}(S_{1}) & \Card \mathcal{K}(T_{S}\cap S_{1}) & K_{1} & 
\Card (\Gamma^{K_{1}}\cap \Gamma^{S_{1}}) \\
\hline
A_{1} & 5 & \emptyset & 1 & 0 & A_{3}\, \hbox{or}\, A_{5} &  1 \\
\hline
A_{2} & 5 & \emptyset & 1 & 0 & A_{5} & 2 \\
\hline
A_{2} & 5 & A_{1} & 1 & 1 & A_{3} & 2 \\
\hline
A_{3} & 6 & A_{1} & 2 & 1 & A_{5} & 3 \\ 
\hline
A_{3} & 5 & A_{1} & 2 & 1 & E_{6} & 3 \\
\hline
A_{4} & 6 & A_{1} & 2 & 1 & E_{6} & 4 \\
\hline
A_{4} & 6 & A_{3} & 2 & 2 & A_{5} & 4 \\ 
\hline 
D_{4} & 5 & A_{3} & 4 & 2 & E_{6} & 4 \\
\hline
D_{5} & 6 & A_{3} & 4 & 2 & E_{6} & 9 \\
\hline
\end{array}
$$
Thus 
$$
\Card \mathcal{K}(S_{1})-\Card \mathcal{K}(T_{S}\cap S_{1})
+ \Card (\Gamma^{S_{1}} \cap \Gamma^{K_{1}} )
\geq 2(\dim V_{S}-\Card \mathcal{K}(\Pi )).
$$
A direct verification gives the result. Note that as in the case of type $A_{\ell}$,
$\mathcal{K}(\Pi )$ is totally ordered by inclusion, so $T_{S}$ is connected.
\end{proof}

\begin{remarks}
1. Theorem \ref{Main} says that if $\mathcal{K}(S) \subset \mathcal{K}(\Pi )$ or equivalently,
$S'\in \mathcal{K}(\Pi )$ for any connected component $S'$ of $S$, 
then $\chi (\mathfrak{p}_{S})+\chi (\mathfrak{u}_{S})=\rank \mathfrak{g}$. 

2. When $\mathfrak{g}$ is of type $B_{\ell}$, $C_{\ell}$, $D_{2n}$, $E_{7}$, $E_{8}$,
$F_{4}$ or $G_{2}$, we have $\Card \mathcal{K}(\Pi )=\rank \mathfrak{g}$. In these cases, 
the condition (i) in Theorem \ref{Main} 
is equivalent to $\mathcal{K}(S)\subset \mathcal{K}(\Pi )$, and consequently, 
condition (ii) is automatically satisfied.
\end{remarks}

\begin{example}
Let us consider the case of minimal parabolic subalgebras. So 
$S=\{ \alpha \}$ and \Eqref{newupperbound} is an equality. It follows that 
$$
\chi (\mathfrak{p}) + \chi (\mathfrak{u} )=\left\{
\begin{array}{ll}
\rank \mathfrak{g} & \hbox{if } \{ \alpha \} \in \mathcal{K}(\Pi ), \\
\rank \mathfrak{g} & \hbox{if } \{ \alpha \} \not \in \mathcal{K}(\Pi ) \hbox{ and }
\dim V_{S}=\Card \mathcal{K}(\Pi )+1, \\
\rank \mathfrak{g} + 2 & \hbox{if } \{ \alpha \} \not \in \mathcal{K}(\Pi ) \hbox{ and }
\dim V_{S}=\Card \mathcal{K}(\Pi ).
\end{array}\right.
$$
Thus the minimal parabolic subalgebras $\mathfrak{p}_{S}$ 
verifying $\chi (\mathfrak{p}_{S}) + \chi (\mathfrak{u}_{S}) = \rank \mathfrak{g}$
are (in the simple roots numbering of \cite[Chapter 18]{TYbook}):
$$
\begin{array}{|c|c|c|}
\hline
\hbox{Type} & \mathcal{K}(S) \not\subset \mathcal{K}(\Pi ) & \mathcal{K}(S)\subset \mathcal{K}(\Pi) \\
\hline
A_{\ell} & \hbox{\vrule height20pt width0pt depth12pt} i\neq {\displaystyle\frac{\ell+1}{2}} &
\hbox{\vrule height20pt width0pt depth12pt} i= {\displaystyle\frac{\ell+1}{2}} \\
\hline
B_{\ell} & \hbox{none} & i \hbox{ odd} \\
\hline
C_{\ell} & \hbox{none} & i=\ell \\
\hline
D_{2n+1} & i=2n, 2n+1 & i < 2n \hbox{ odd} \\
\hline
D_{2n} & \hbox{none} & i \hbox{ odd or } i=2n \\
\hline
E_{6} & i\neq 2,4 &  i=4 \\
\hline
E_{7} & \hbox{none} & i=2,3,5,7 \\
\hline
E_{8} & \hbox{none} & i=2,3,5,7 \\
\hline
F_{4} & \hbox{none} & i=2 \\
\hline
G_{2} & \hbox{none} & i=1  \\
\hline
\end{array}
$$
\end{example}

\begin{example}
In the other extremity, it is easy to check that maximal parabolic subalgebras of 
$\mathfrak{g}=\mathrm{sl}_{\ell+1}$ verifying 
$\chi (\mathfrak{p})+\chi (\mathfrak{u})=\rank \mathfrak{g}$
are exactly the ones associated to simple roots at the extremities 
of the Dynkin diagram.
\end{example}

\end{document}